\RequirePackage{ifpdf}
\ifpdf 
\documentclass[pdftex]{sigma}
\else
\documentclass{sigma}
\fi

\begin{document}
\allowdisplaybreaks

\renewcommand{\PaperNumber}{034}

\FirstPageHeading

\ShortArticleName{On Orthogonality Relations for Dual Discrete
$q$-Ultraspherical Polynomials}

\ArticleName{On Orthogonality Relations\\
 for Dual Discrete $\boldsymbol{q}$-Ultraspherical Polynomials}

\Author{Valentyna A. GROZA~$^\dag$ and Ivan I. KACHURYK~$^\ddag$}

\AuthorNameForHeading{V.A. Groza and I.I. Kachuryk}

\Address{$^\dag$~National Aviation University, 1 Komarov Ave.,
Kyiv,  03058 Ukraine}
\EmailD{\href{mailto:groza@i.com.ua}{groza@i.com.ua}}

\Address{$^\ddag$~Khmel'nyts'kyi National University,
Khmel'nyts'kyi, Ukraine}
\EmailD{\href{mailto:kachuryk@ief.tup.km.ua}{kachuryk@ief.tup.km.ua}}

\ArticleDates{Received February 14, 2006, in f\/inal form February
28, 2006; Published online March 16, 2006}

\bigskip

\leftline{\small \bfseries \itshape Submitted by Anatoly Klimyk}

\Abstract{The dual discrete $q$-ultraspherical polynomials
$D_n^{(s)}(\mu (x;s)|q)$ correspond to indeterminate moment
problem and, therefore, have one-parameter family of extremal
ortho\-go\-nality relations. It is shown that special cases of dual
discrete $q$-ultraspherical polynomials $D_n^{(s)}(\mu (x;s)|q)$,
when $s=q^{-1}$ and $s=q$, are directly connected with
$q^{-1}$-Hermite polynomials. These connections are given in an
explicit form. Using these relations, all extremal orthogonality
relations for these special cases of polynomials $D_n^{(s)}(\mu
(x;s)|q)$ are found.}


\Classification{33D45; 81Q99}

\section{Introduction}

Quantum groups and $q$-deformed variants of the quantum harmonic
oscillator were introduced approximately two decades ago.
Representations of quantum groups and realizations of
$q$-oscillators are closely related to basic hypergeometric
functions and $q$-orthogonal polynomials. Instead of Jacobi,
Gegenbauer (ultraspherical), Hermite polynomials, in the theory of
quantum groups and under studying $q$-oscillators we deal with
dif\/ferent kinds of $q$-Jacobi polynomials, of $q$-ultraspherical
polynomials and of $q$-Hermite polynomials.

It is known that $q$-orthogonal polynomials are more complicated
than classical orthogonal polynomials. It is interesting that main
development of the theory of $q$-orthogonal polynomials in fact
coincides in time with appearing of quantum groups (see, for
example, \cite{1} and \cite{2}). Due to great importance for
contemporary mathematical and theoretical physics the theory of
$q$-orthogonal polynomials is under extensive development. A great
attention is paid to those of them, which are directly related to
representations of quantum groups and to realizations of
$q$-oscillators.

In this paper we deal with the so-called $q^{-1}$-Hermite
polynomials (which are closely related to the
Biedenharn--Macfarlane oscillator; see \cite{3}) and with the dual
discrete $q$-ultraspherical polynomials. The $q^{-1}$-Hermite
orthogonal polynomials were discovered by Askey \cite{4} and the
discrete $q$-ultraspherical polynomials and their duals were
introduced in \cite{5}.

The $q^{-1}$-Hermite orthogonal polynomials correspond to
indeterminate moment problem. This means that they have a
one-parameter family of orthogonality relations (these
orthogonalities are important for studying the
Biedenharn--Macfarlane oscillator when $q>1$; see \cite{3}). These
orthogonality relations were explicitly derived in \cite{6}. It
was shown in \cite{5} that dual $q$-ultraspherical polynomials
also correspond to indeterminate moment problem and, therefore,
there exists a~one-parameter family of orthogonality relations for
them. In \cite{5}, there were derived only two ortho\-gonality
relations for these polynomials. Other orthogonality relations for
them are not known. As in the case of usual ultraspherical
(Gegenbauer) polynomials, these orthogonality relations are very
important.

Let us remark that f\/inding extremal orthogonality relations for
orthogonal polynomials, which correspond to indeterminate moment
problems, is a very dif\/f\/icult task. In fact, at present  only one
example of such measures is known; that is extremal orthogonality
relations for the $q^{-1}$-Hermite polynomials.

In this paper we derive explicit formulas that connect the dual
$q$-ultraspherical polynomials $D_n^{(s)}(\mu (x;s)|q)$, when
$s=q^{-1}$ and $s=q$, with the $q^{-1}$-Hermite polynomials, and
this gives a possibility to obtain all extremal measures for these
particular types of the the dual $q$-ultraspherical polynomials.

Note that the dual $q$-ultraspherical polynomials $D_n^{(s)}(\mu
(x;s)|q)$ are related to the discrete series representations of
the quantum Drinfeld--Jimbo algebra $U_q({\rm su}_{1,1})$;
see~\cite{7}.

Below we use (without additional explanation) notations of the
theory of $q$-special functions (see \cite{8}). In particular, the
shifted $q$-factorials $(a;q)_n$ are def\/ined as
\[
(a;q)_n=(1-a)(1-aq)\big(1-aq^2\big)\cdots \big(1-aq^{n-1}\big).
 \]
The basic hypergeometric functions are given by the formula
\[
{}_r\phi_{r-1}\left( a_1, a_2,\ldots, a_r; b_1, b_2,\ldots, b_s;
q,x \right)=\sum_{n=0}^\infty
 \frac{(a_1;q)_n (a_2;q)_n\ldots (a_r;q)_n}
 {(b_1;q)_n (b_2;q)_n\ldots (b_{r-1};q)_n} \frac{x^n}{(q;q)_n}.
 \]

\section[Dual discrete $q$-ultraspherical polynomials
and $q^{-1}$-Hermite polynomials]{Dual discrete $\boldsymbol{q}$-ultraspherical polynomials\\
and $\boldsymbol{q^{-1}}$-Hermite polynomials}

The discrete $q$-ultraspherical polynomials $C_n^{(s)}(x;q)$ are
the polynomials
\[
C_n^{(s)}(x;q)= {}_3\phi_2\left( q^{-n}, -sq^{n+1}, x; \sqrt{s}q,
-\sqrt{s}q;  q,q \right).
\]
These polynomials are orthogonal on a countable set of points (see
formula~(4) in~\cite{5}).

In \cite{9}, the procedure of creation of dual sets of
polynomials, orthogonal on a discrete set of points, is being
developed. The polynomials $D_n^{(s)}(\mu(x;s);q)$ dual to the
discrete $q$-ultraspherical polynomials $C_n^{(s)}(x;q)$ are of
particular case of the dual big $J$-Jacobi polynomials (see
\cite{10}). They are given by the formula
\[
D_n^{(s)}(\mu(x;s);q)= {}_3\phi_2\left( q^{-x}, sq^{x+1}, q^{-n};
\sqrt{s}q, -\sqrt{s}q;  q, -q^{n+1} \right),
\]
where $\mu(x;s)=q^{-x}+sq^{x+1}$. They are polynomials in
$\mu(x;s)$. These polynomials satisfy the three-term recurrence
relation
 \begin{gather}
\big(q^{-x}+sq^{x+1}\big) D_n^{(s)}(\mu (x;s)|q)=
-q^{-2n-1}\big(1-sq^{2n+2}\big) D_{n+1}^{(s)}(\mu (x;s)|q)
 \nonumber\\
   \label{rec}
\qquad {}  +q^{-2n-1}(1+q)D_n^{(s)}(\mu (x;s)|q)
-q^{-2n}\big(1-q^{2n}\big)D_{n-1}^{(s)}(\mu (x;s)|q).
  \end{gather}
They correspond to indeterminate moment problem and it means that
they satisfy a one-parameter family of orthogonality relations.
Only two orthogonality relations were found in~\cite{5}:
 \begin{gather}
\sum_{k=0}^\infty
\frac{(1{-}sq^{4k+1})(sq;q)_{2k}}{(1{-}sq)(q;q)_{2k}}\,
q^{k(2k-1)}
D_n^{(s)}(\mu (2k)|q) D_{n'}^{(s)}(\mu (2k)|q)\nonumber\\
\qquad {} =\frac{(sq^3;q^2)_\infty}{(q;q^2)_\infty}
 \frac{(q^2;q^2)_nq^{-n}}{(sq^2;q^2)_n}\delta_{nn'} ,\label{ort-1}
 \\
 \sum_{k=0}^\infty
\frac{(1{-}sq^{4k+3})(sq;q)_{2k+1}}{(1{-}sq)(q;q)_{2k+1}}\,
q^{k(2k+1)} D_n^{(s)}(\mu (2k+1)|q) D_{n'}^{(s)}(\mu (2k+1)|q)\nonumber\\
  \label{ort-2}
\qquad{}  =\frac{(sq^3;q^2)_\infty}{(q;q^2)_\infty}
 \frac{(q^2;q^2)_n q^{-n}}{(sq^2;q^2)_n} \delta_{nn'},
  \end{gather}
where $\mu(2k)\equiv \mu(2k;s)$, $\mu(2k+1)\equiv \mu(2k+1;s)$ and
$0<s<q^{-2}$.

We shall need also  $q^{-1}$-Hermite polynomials $h_n(x|q)$, which
are given as
 \[
h_n(\sinh \varphi| q)=\sum_{k=0}^n \frac{(-1)^k q^{k(k-n)} (q;
q)_n}{( q; q)_k( q; q)_{n-k}}e^{\varphi(n-2k)}.
 \]
They are def\/ined as $h_n(x|q)={\rm i}^{-n}H_n({\rm i}x|q^{-1})$,
where $H_n(y|q)$ are well-known continuous $q$-Hermite
polynomials. The polynomials $h_n(x|q)$ satisfy the following
recurrence relation
  \begin{gather}  \label{rec-3}
h_{n+1}(x|{ q})+ {q}^{-n}\big(1-{q}^n\big) h_{n-1}(x|{q}))=2x
h_{n}(x|{q}).
  \end{gather}
It follows from this relation that the polynomials $h_{2k}(x|q)$
are even:
\[
h_{2k}(-x|q)=h_{2k}(x|q)
\]
and the polynomials $h_{2k+1}(x|q)$ have the form
\[
h_{2k+1}(x|q)=x\tilde h_{2k}(x|q) ,
\]
where $\tilde h_{2k}(x|q)$ are even polynomials: $\tilde
h_{2k}(-x|q)=\tilde h_{2k}(x|q)$.

The $q^{-1}$-Hermite polynomials $h_n(x|q)$ correspond to
indeterminate moment problem and, therefore, have one-parameter
family of extremal orthogonality measures, which were given
in~\cite{6}. They are given by the parameter $a$, $q\le a<1$, and
have the form
  \begin{gather}  \label{H-1}
 \sum_{m=-\infty}^\infty
\frac{a^{4m}q^{m(2m-1)}(1+a^2q^{2m})}{(-a^2;q)_\infty
(-q/a^2;q)_\infty (q;q)_\infty}
h_n(x_m|q)h_{n'}(x_m|q)=q^{-n(n+1)/2}{(q;q)_n} \delta_{nn'},
  \end{gather}
where
 \[
x_m=\frac12 \big(a^{-1}q^{-m}-aq^m\big).
\]
The relations \eqref{H-1} were proved in~\cite{6} only for $a$
from the interval $q<a<1$. It is only said there that the
relation~\eqref{H-1} is true also for $a=q$. Below we shall give a
proof of~\eqref{H-1} for $a=q$.

\medskip

\noindent {\bf Proposition.} {\it The following expression for
$q^{-1}$-Hermite polynomials $h_n(x|q)$ hold:
  \begin{gather}  \label{rel-1}
 h_{2k}(\sinh
\varphi|q)=(-1)^kq^{-k^2}\big(q;q^2\big)_k D^{(q^{-1})}_k
\big(e^{2\varphi}+e^{-2\varphi}|q\big),
  \\
    \label{rel-2}
h_{2k+1}(\sinh
\varphi|q)=(-1)^kq^{-k(k+1)}\big(q^3;q^2\big)_k(2\sinh \varphi)
D^{(q)}_k \big(qe^{2\varphi}+qe^{-2\varphi}|q\big),
  \end{gather}
 where $k$ are
nonnegative integers and $D^{(s)}_k (\mu(m;s)|q)$ are dual
discrete $q$-ultraspherical polynomials.}

 \medskip

 \begin{proof} We take the recurrence relation \eqref{rec-3} for
$n=2k$:
\[
2x h_{2k}(x|{q})= h_{2k+1}(x|{ q})+ {q}^{-2k}(1-{q}^{2k})
h_{2k-1}(x|{q}))
\]
and substitute there the expressions for $h_{2k+1}(x|{ q})$ and
$h_{2k-1}(x|{ q})$, obtained from \eqref{rec-3}. As a~result, we
get the relation
\begin{gather*}
(2x)^2 h_{2k}(x|{q})= h_{2k+2}(x|{q})+\big[q^{-2k}\big(1+q^{-1}\big)-2\big]h_{2k}(x|{q})\\
\phantom{(2x)^2 h_{2k}(x|{q})=}{} +
{q}^{-4k+1}\big(1-{q}^{2k}\big) \big(1-{q}^{2k-1}\big)
h_{2k-2}(x|{q})
\end{gather*}
which can be written in the form
\begin{gather}
\big(e^{2\varphi}+e^{-2\varphi}\big) h_{2k}(\sinh \varphi|{q})=
h_{2k+2}(\sinh \varphi|{q})+q^{-2k}\big(1+q^{-1}\big)h_{2k}(\sinh \varphi|{q})\nonumber\\
  \label{Herm}
 \phantom{\big(e^{2\varphi}+e^{-2\varphi}\big) h_{2k}(\sinh \varphi|{q})=}{}
 +{q}^{-4k+1}\big(1-{q}^{2k}\big)
\big(1-{q}^{2k-1}\big) h_{2k-2}(\sinh \varphi|{q}) .
  \end{gather}

Now we take the recurrence relation \eqref{rec} for $s=q^{-1}$:
\begin{gather*}
y
D^{(q^{-1})}_n(y|q)=-q^{-2n-1}\big(1-q^{2n+1}\big)D^{(q^{-1})}_{n+1}(y|q)
+ q^{-2n-1}(1+q) D^{(q^{-1})}_n(y|q)
 \\
\phantom{y D^{(q^{-1})}_n(y|q)=}{} -
q^{-2n}\big(1-q^{2n}\big)D^{(q^{-1})}_{n-1}(y|q)
\end{gather*}
and substitute there the expression
\[
D^{(q^{-1})}_n(y|q)=(-1)^nq^{n^2}\big(q;q^2\big)^{-1}_n \tilde
D^{(q^{-1})}_n(y|q)
\]
(compare it with the f\/irst relation in the proposition), where
$\tilde D^{(s)}_n(y|q)$ are polynomials multiple to the
polynomials $D^{(s)}_n(y|q)$. As a result, we obtain the relation
 \begin{gather}
y  \tilde D^{(q^{-1})}_n(y|q)=\tilde D^{(q^{-1})}_{n+1}(y|q) +
q^{-2n-1}(1+q) \tilde D^{(q^{-1})}_n(y|q)
 \nonumber\\
   \label{dual}    \phantom{y  \tilde D^{(q^{-1})}_n(y|q)=}{}
 +q^{-4n+1}\big(1-q^{2n}\big)\big(1-q^{2n-1}\big) \tilde D^{(q^{-1})}_{n-1}(y|q).
  \end{gather}
Setting here $y=e^{2\varphi}+e^{-2\varphi}$ and then comparing the
relations~\eqref{Herm} and~\eqref{dual} we get
 \begin{gather}  \label{equ}
h_{2k}(\sinh \varphi|q)=c\tilde D^{(q^{-1})}_k
\big(e^{2\varphi}+e^{-2\varphi}|q\big),
 \end{gather}
where $c$ is a constant. In order to f\/ind this constant we put
$k=0$ in \eqref{equ}. Since $h_0(x|q)=1$ and $\tilde
D^{(q^{-1})}_0(y|q)=1$, we obtain that $c=1$ and the f\/irst
relation of the proposition is true.

The relation \eqref{rel-2} is proved similarly. We take the
recurrence relation \eqref{rec-3} for $n=2k+1$:
\[
2x h_{2k+1}(x|{q})= h_{2k+2}(x|{ q})+
{q}^{-2k-1}\big(1-{q}^{2k+1}\big) h_{2k}(x|{q})
\]
and substitute there the expressions for $h_{2k+2}(x|{ q})$ and
$h_{2k}(x|{ q})$, obtained from~\eqref{rec-3}. As a result, we get
the relation which can be written in the form
\begin{gather}
\big(qe^{2\varphi}+qe^{-2\varphi}\big) h_{2k+}(\sinh \varphi|{q})=
q
h_{2k+3}(\sinh \varphi|{q})\nonumber\\
  \label{Herm-2}
 \qquad{}
+q^{-2k}\big(1+q^{-1}\big)h_{2k+1}(\sinh \varphi|{q})
+{q}^{-4k}\big(1-{q}^{2k}\big) \big(1-{q}^{2k+1}\big) h_{2k}(\sinh
\varphi|{q}).
  \end{gather}

Now we take the recurrence relation \eqref{rec} for $s=q$:
\begin{gather*}
y D^{(q)}_n(y|q)=-q^{-2n-1}\big(1-q^{2n+3}\big)D^{(q)}_{n+1}(x|q)
+ q^{-2n-1}(1+q) D^{(q)}_n(y|q)
 \\
\phantom{y D^{(q)}_n(y|q)=}{}  -
q^{-2n}\big(1-q^{2n}\big)D^{(q)}_{n-1}(y|q)
\end{gather*}
and substitute there the expression
\[
D^{(q)}_n(y|q)=(-1)^nq^{n(+1)}\big(q^3;q^2\big)^{-1}_n \tilde
D^{(q)}_n(y|q)
\]
(compare with the second relation in the proposition). As a
result, we obtain the relation
\begin{gather}
y  \tilde D^{(q^{-1})}_n(y|q)=q \tilde D^{(q^{-1})}_{n+1}(y|q) +
q^{-2n-1}(1+q) \tilde D^{(q^{-1})}_n(y|q)\nonumber\\
  \label{dual-2}\phantom{y  \tilde D^{(q^{-1})}_n(y|q)=}{}
 +q^{-4n}\big(1-q^{2n}\big)(1-q^{2n+1}) \tilde
D^{(q^{-1})}_{n-1}(y|q).
  \end{gather}
Setting here $y=qe^{2\varphi}+qe^{-2\varphi}$ and then comparing
the relations~\eqref{Herm-2} and~\eqref{dual-2} we get
 \begin{gather}  \label{equ-2}
h_{2k+1}(\sinh \varphi|q)=c\tilde D^{(q)}_k
\big(e^{2\varphi}+e^{-2\varphi}|q\big),
 \end{gather}
where $c$ is a constant. Putting $k=0$ in \eqref{equ-2}, we f\/ind
that $c=1$ and the second relation of the proposition is true. The
proposition is proved.
\end{proof}

\section[Orthogonality relations for $h_n(x|q)$ when $a=q$]{Orthogonality relations
for $\boldsymbol{h_n(x|q)}$ when $\boldsymbol{a=q}$}

The aim of this section is to show that the orthogonality
relations \eqref{ort-1} and \eqref{ort-2} are equivalent to to the
orthogonality relation \eqref{H-1} for $a=q$. Thus, we give a
proof of the orthogonality relation~\eqref{H-1} for $a=q$,
starting from the orthogonality relations~\eqref{ort-1}
and~\eqref{ort-2}.

It follows from relation \eqref{ort-1} that the dual discrete
$q$-ultraspherical polynomials $D^{(q^{-1})}_n(x|q)$ satisfy the
orthogonality relation
 \begin{gather}  \label{ort-3}
\sum_{m=0}^\infty
\big(1+q^{2m}\big)q^{m(2m-1)}D_n^{(q^{-1})}(y_m|q)
D_{n'}^{(q^{-1})}(y_m|q)=\frac{(q^2;q^2)_\infty}{(q;q^2)_\infty}
\frac{(q^2;q^2)_nq^{-n}}{(q;q^2)_n} \delta_{nn'},
 \end{gather}
where $y_m=q^{2m}{+}q^{-2m}$. By using Proposition, we express
$D^{(q^{-1})}_n(x_m|q)$ in terms of $h_{2n}(\frac12
(q^{-m}-q^{m})|q)$ and substitute into \eqref{ort-3}. As a result,
we obtain
 \begin{gather}  \label{ort-4}
\sum_{m=0}^\infty \big(1+q^{2m}\big)q^{m(2m-1)}h_{2n}(x_m|q)
h_{2n'}(x_m|q)=\frac{(q^2;q^2)_\infty}{(q;q^2)_\infty}
q^{-n(2n+1)}(q;q)_{2n} \delta_{nn'},
 \end{gather}
where $x_m=\frac12 (q^{-m}-q^{m})$. Note that
\begin{gather}
 \frac{(q^2;q^2)_\infty}{(q;q^2)_\infty}=(-q;q)^2_\infty (q;q)_\infty
 =\frac12 (-1;q)_\infty (-q;q)_\infty (q;q)_\infty
 \nonumber\\
  \label{rel-q}\phantom{\frac{(q^2;q^2)_\infty}{(q;q^2)_\infty}}{}
  =2q^{-1} \big(-q^2;q\big)_\infty \big(-q^{-1};q\big)_\infty
(q;q)_\infty.
 \end{gather}

The expression $(1+q^{2m})q^{m(2m-1)}$ on the right hand side of
\eqref{ort-4} does not change at replacing~$m$ by $-m$. Since
$h_{2n}(x|q)$ is an even function, then due to \eqref{rel-q} the
relation \eqref{ort-4} can be written~as
\begin{gather}  \label{ort-5}
\sum_{m=-\infty}^\infty
\frac{(1+q^{2m})q^{-1}q^{m(2m-1)}}{(-q^2;q)_\infty
(-q^{-1};q)_\infty (q;q)_\infty}
 h_{2n}(x_m|q) h_{2n'}(x_m|q)= q^{-n(2n+1)}(q;q)_{2n}
 \delta_{nn'}.
 \end{gather}
This relation is equivalent to the orthogonality relation
\eqref{H-1} at $a=q$ if $n$ and $n'$ are replaced there by $2n$
and $2n'$, and $m$ by $-m-1$, respectively.

Now we use the relation \eqref{ort-2} for the dual discrete
$q$-ultraspherical polynomials $D^{(q)}_n(x|q)$. Then after
shifting $k$ to $k-1$ we have
\begin{gather}
\sum_{m=1}^\infty
\big(1+q^{2m}\big)\big(1-q^{2m}\big)^2q^{(m-1)(2m-1)}D_n^{(q)}(y_m|q)
D_{n'}^{(q)}(y_m|q)
\nonumber\\
 \label{ort-6}    \qquad{}
=\frac{(q^2;q^2)_\infty}{(q;q^2)_\infty}
\frac{(1-q)^2(q^2;q^2)_nq^{-n}}{(1-q^{2n+1})(q;q^2)_n}
\delta_{nn'},
 \end{gather}
where $y_m=q^{2m+1}+q^{-2m+1}$. By using the above Proposition, we
express $D^{(q)}_n(y_m|q)$ in terms of $h_{2n+1}(\frac12
(q^m-q^{-m})|q)$ and substitute into \eqref{ort-6}. We obtain
 \begin{gather}
\sum_{m=1}^\infty \big(1+q^{2m}\big)q^{m(2m-1)}h_{2n+1}(x_m|q)
h_{2n'+1}(x_m|q)\nonumber\\
\qquad {}=\frac{(q^2;q^2)_\infty}{(q;q^2)_\infty}
q^{-(n+1)(2n+1)}(q;q)_{2n+1}\delta_{nn'}, \label{ort-7}
 \end{gather}
where, as before, $x_m=\frac12 (q^{-m}-q^{m})$. Since
$h_{2n+1}(0|q)=0$, then the summation $\sum\limits_{m=1}^\infty$
here can be replaced by $\sum\limits_{m=0}^\infty$, and then by
$\frac12 \sum\limits_{m=-\infty}^\infty$ (since $h_{2n+1}(x|q)$ is
an odd function). As in the case of \eqref{ort-4}, it is easy to
show that the last relation is equivalent to the relation
\eqref{H-1} if we replace $m$ by $m-1$, and $n$ and $n'$ by $2n+1$
and $2n'+1$, respectively.

Since the polynomials $h_{2n}(x_m|q)$ are even and the polynomials
$h_{2n+1}(x_m|q)$ are odd, then
\begin{gather*}
\sum_{m=0}^\infty \big(1+q^{2m}\big)q^{m(2m-1)}h_{2n}(x_m|q)
h_{2n'+1}(x_m|q)
 \\ \qquad{}=- \sum_{m=-\infty}^0 \big(1+q^{2m}\big)q^{m(2m-1)}h_{2n}(x_m|q)
h_{2n'+1}(x_m|q).
\end{gather*}
Therefore,
 \begin{gather}  \label{ort-8}
\sum_{m=-\infty}^\infty \big(1+q^{2m}\big)q^{m(2m-1)}h_{2n}(x_m|q)
h_{2n'+1}(x_m|q)=0.
 \end{gather}
This gives the mutual orthogonality of the polynomials
$h_{2n}(x_m|q)$, $n=0,1,2,\ldots$, to the polynomials
$h_{2n+1}(x_m|q)$, $n=0,1,2,\ldots$. Thus, we proved the
relation~\eqref{H-1} for $a=q$.

\section[Orthogonality relations for the polynomials $D^{(q^{-1})}_n$
and $D^{(q)}_n$]{Orthogonality relations for the polynomials
$\boldsymbol{D^{(q^{-1})}_n}$ and $\boldsymbol{D^{(q)}_n}$}

In this section we use the orthogonality relations for the
$q^{-1}$-Hermite polynomials in order to f\/ind orthogonality
relations for the dual discrete $q$-ultraspherical polynomials
$D^{(q^{-1})}_n$ and $D^{(q)}_n$. We substitute the expression
\eqref{rel-1} for $h_{2k}(\sinh \varphi|q)$ into the orthogonality
relation \eqref{H-1}. After some transformations we receive the
orthogonality relations for $D^{(q^{-1})}_n$:
 \begin{gather}
 \sum_{m=-\infty}^\infty a^{4m+1}q^{2m^2}\big(a^{-1}q^{-m}+aq^{m}\big)
D^{(q^{-1})}_k(y_m|q) D^{(q^{-1})}_k(y_m|q)\nonumber\\
\qquad{} =\frac{(-a^2;q)_\infty (-q/a;q)_\infty
(q;q)_\infty}{q^{k}(q;q)^{-1}_{2k}(q;q^2)^{2}_k}
\delta_{kk'},\label{ort-9}
 \end{gather}
where
\[
y_m=a^{-2}q^{-2m}+a^2q^{2m}.
\]
These relations are parametrized by $a$, $q\le a<1$. For dif\/ferent
values of $a$ the corresponding orthogonality relations are not
equivalent.

Now we substitute into \eqref{H-1} the relation \eqref{rel-2} for
$h_{2k+1}(\sinh \varphi|q)$. After some transformation we obtain
 \begin{gather}
 \sum_{m=-\infty}^\infty a^{4m}q^{m(2m-1)}\big(a^{-2}q^{-2m}-a^2q^{2m}\big)
D^{(q)}_k(y_m|q) D^{(q)}_k(y_m|q) \nonumber\\
\qquad{}=\frac{(-a^2;q)_\infty (-q/a;q)_\infty
(q;q)_\infty}{q^{k+1}(q;q)^{-1}_{2k+1}(q^3;q^2)^{2}_k}
\delta_{kk'}\label{ort-10}
 \end{gather}
where
 \[
y_m=\big(a^{-2}q^{-2m}+a^2q^{2m}\big)q.
 \]
These relations are also parametrized by $a$, $q\le a<1$.
Orthogonality relations, corresponding to dif\/ferent values of $a$,
are not equivalent.

The orthogonality relations \eqref{H-1} for the $q^{-1}$-Hermite
polynomials are extremal, that is, the set of polynomials
$h_n(x_m|q)$, $n=0,1,2,\ldots$, is complete in the corresponding
Hilbert space (see~\cite{11}). Therefore, the orthogonality
relations \eqref{ort-9} and \eqref{ort-10} for the polynomials
$D^{(q^{-1})}_k$ and $D^{(q)}_k$, respectively, are extremal.
Moreover, since the set of orthogonality relations~\eqref{H-1} is
a complete set of extremal orthogonality relations for the
$q^{-1}$-Hermite polynomials, then the set of
relations~\eqref{ort-9} is a~complete set of extremal
orthogonality relations for the polynomials~$D^{(q^{-1})}_k$ and
the set of relations \eqref{ort-10} is a complete set of extremal
orthogonality relations for the polynomials~$D^{(q)}_k$.

\section{Concluding remarks}

The dual discrete $q$-ultraspherical polynomials $D^{(s)}_n(\mu
(x;s)|q)$ correspond to the indeterminate moment problem, that is,
there are inf\/initely many (one-parameter family) of orthogonality
relations for them. We have found these orthogonality relations
for the cases when $s=q$ and $s=q^{-1}$. Namely, we have reduced
these cases to the $q^{-1}$-Hermite polynomials for which
orthogonality relations are known. It leaves  unsolved the problem
of f\/inding orthogonality relations for other cases. We believe
that the orthogonality relations~\eqref{ort-9} and~\eqref{ort-10}
will help in f\/inding orthogonality relations for these other
cases.

There exists another class of dual discrete $q$-ultraspherical
polynomials (they are denoted as  $\tilde D^{(s)}_n(\mu
(x;-s)|q)$; see \cite[Section~1]{5}). They also correspond to the
indeterminate moment problem. It is interesting to f\/ind
orthogonality relations for them. However, they cannot be reduced
to the $q^{-1}$-Hermite polynomials. It is shown in~\cite{5} that
they are reduced to Berg--Ismail polynomials. This connection
gives an inf\/inite number of orthogonality relations for them.
However, they are not extremal orthogonality relations. So, the
problem of f\/inding these orthogonality relations are unsolved.

It is well-known (see \cite{3}) that $q^{-1}$-Hermite polynomials
are closely related to the Biedenharn--Macfarlane $q$-oscillator.
Thus, the polynomials, considered in this paper, is also related
to this $q$-oscillator. Then the following question appears: Are
the polynomials $\tilde D^{(s)}_n(\mu (x;-s)|q)$ related to the
$q$-oscillator? Probably, they correspond to some simple
deformation of the $q$-oscillator.

\LastPageEnding

\end{document}